# Three classes of Ermakov systems and nonlocal symmetries

**ARUNAYE, F.I.**

*Correspondence address*: Department of mathematics and Computer Science, The University of the West IndiesMona, Kingston 7. Jamaica.

*Permanent address*: Department of Mathematics and computer Science Delta State University Abraka. Nigeria

*E-mail*: **ifarunaye@yahoo.com**; **festus.arunaye@uwimona.edu.jm**

## Abstract

Ermakov systems have attracted enormous treatments in recent times particularly in symmetry analysis. In this paper we consider three classes of the Ermakov systems by using a simple algebraic reduction process with imposed conditions on the magnitude of the angular momentum of each system class to obtain new generalized symmetries. We note that this imposed condition transforms the Kepler-Ermakov systems to the generalized Ermakov systems.

## 1.0     Introduction

The possession of first integral and three Lie point symmetry generators of the algebra $sl(2,R)$ that are characteristics of Ermakov systems of second-order ordinary differential equations are well known in the literature (Ermakov [1], Harin [2], Athorne [3], Goedert and Haas [4], Haas and Goedert [5], Simic [6], Leach and Karasu (Kalkanli) [7], Goodall and Leach [8]). The absence of sufficient number of Lie point symmetries for the Kepler problem in the context of complete symmetry group of dynamical systems brought to the fore the introduction of the nonlocal symmetries of dynamical systems by Krause [9]. Nucci [10] introduced her concept of reduction order combined with the Lie algorithm for obtaining the classical symmetries of differential equations to obtain the complete symmetry group of the Kepler problem according to Krause [9] and as well voided the earlier assertions of Krause [9] that these nonlocal symmetries could not be obtained by Lie algorithm. The Nucci [10] reduction process became so famous in the literature hitherto Arunaye [11] announced a simpler reduction process for reducing dynamical systems to systems of oscillator and equation of motion which admits Lie algorithm for the computation of their infinitesimal vector fields, meanwhile it is already established in Leach, Andriopoulos and Nucci [12] that the Ermanno-Bernoulli constants of dynamical systems are most suitable for reducing dynamical systems to systems of oscillator and equation of motion. Arunaye and White [13] reported alternative constants which were constructed from the Hamilton vector of dynamical systems that are equivalent to the Ermanno-Bernoulli constants in two dimensions and which seemed to provide less cumbersome reduction variables in three-dimensions than the Ermanno-Bernoulli constants Arunaye [14].

The central feature of the Ermakov systems is their property of always having first integrals Haas and Goedert [5], Simic [6] and that this invariant plays a central role in the linearization of Ermakov systems Ray and Reid [15], Haas and Goedert [5], Athorne [3]. The Kepler-Ermakov systems referred to the perturbations of the classical Kepler problem or an autonomous Ermakov system was investigated by Karasu (Kalkanli) and Yildirim [16] and found that these systems are the usual Ermakov systems with frequency function depending on the dynamical variables. Leach and Karasu (Kalkanli) [7] supplemented the analysis of Karasu (Kalkanli) and Yildirim [16] by correcting some results as well as carried out investigation on the same dynamics in which an equivalent transformation of the Kepler-Ermakov systems to new time and rescaled radial distance so that in the discussion of the Kepler-Ermakov systems it suffices to study its polar equivalent system. The paper of Leach, Karasu (Kalkanli), Nucci and Andriopoulos [17] considered the Ermakov's superintegrable-toy for its nonlocal symmetries and asserted the insufficient Lie point symmetries and the unstable algebra $sl(2,R)$ for the complete specification of the system. Also the method of Nucci reduction process brought the representation of the complete symmetry group to the fore, four of which are nonlocal symmetries and the algebra is the direct sum of a one-dimensional Abelian algebra and the semidirect sum of a solvable algebra with a two-dimensional Abelian algebra $[A_1 \oplus \{A_2 \oplus_s 2A_1\}]$ in the notation of the Mubarrakzyanov classification scheme Morozov [18], Mubarakzyanov [19, 20, 21]. Arunaye [22] considered the reduction of three classes of the Ermakov systems by the method of Arunaye [11] and obtained some new nonlocal symmetry for these classes of Ermakov systems. In this paper we utilized the results of Arunaye [22] and imposed certain conditions on the magnitude of the angular momentum (which is not constant in these circumstances) and with a specific function $H$ as a Kernel transformation on the radial component of motion of the Kepler-Ermakov systems. The paper is organized as following. Section 2 recalled the reduction of the three classes of the Ermakov systems to systems of linear second order and conservation law. Section 3 is devoted to the symmetry analysis and finally section 4 is conclusions.

## 1.1 Classes of Ermakov systems

The three classes of the Ermakov systems under consideration are given by

$$\ddot{x} + w^2(t)x = -\frac{x}{r^3}H + \frac{1}{x^3}f(\frac{y}{x})$$
$$\ddot{y} + w^2(t)y = -\frac{y}{r^3}H + \frac{1}{y^3}g(\frac{y}{x}); \qquad (1.1)$$

$$\ddot{x} + w^2(t)x = \frac{1}{yx^2}f(\frac{y}{x})$$
$$\ddot{y} + w^2(t)y = \frac{1}{xy^2}g(\frac{y}{x}) ; \qquad (1.2)$$

and

$$\ddot{x} + w^2(t)x = \frac{1}{x^3}$$
$$\ddot{y} + w^2(t)y = \frac{1}{y^3}, \qquad (1.3)$$

where $f$ and $g$ are arbitrary functions of their arguments and $H$ is a function of unspecified form of dependence upon $x$, $y$ and $r$, denoted specifically by $H = \frac{1}{4}Cr^3 - \frac{1}{r\cos\theta}h(\cot\theta)$ where $C$ is an arbitrary constant (Leach and Karasu (Kalkanli) [17]). Systems (1.1), (1.2) and (1.3) are known as Kepler-Ermakov systems, generalized Ermakov systems and the Ermakov-toy systems respectively. In the symmetry analysis parlance, Leach and Karasu (Kalkanli) [7], Leach, Karasu (Kalkanli), Nucci and Andriopoulos [17] established the plane polar coordinates of (1.1.), (1.2) and (1.3) for their generalized symmetry analysis where $x = r\cos\theta$ and $y = r\sin\theta$.

## 2.0 On the reduction of classes of Ermakov systems

In the following subsections we highlight the reduction process for reducing dynamical systems to systems of second order linear equation and an equation of motion. The method of Arunaye [11] is applied to these three classes of Ermakov systems (1.1), (1.2) and (1.3) in order to utilize the Lie point symmetry analysis method to these classes of dynamical systems and obtain their generalized symmetries.

## 2.1 Reduction process of the Kepler-Ermakov systems

The so called Kepler-Ermakov systems studied by Leach and Karasu (Kalkanli) [7] in polar system have the radial and transversal components of the motion respectively given by

$$\ddot{r} - r\dot{\theta}^2 = \frac{1}{r^3\cos\theta}h(\cot\theta) + \frac{1}{r^3}\{\sec^2\theta f(\tan\theta) + \cosec^2\theta g(\tan\theta)\}, \qquad (2.1)$$

$$r\ddot{\theta} + 2\dot{r}\dot{\theta} = -\frac{1}{r^3}\{\sec^2\theta\tan\theta f(\tan\theta) - \cos ec^2\theta\cot\theta g(\tan\theta)\}. \tag{2.2}$$

Now from (2.2) we have that (Arunaye [22])

$$(r^4\dot{\theta}^2)^{\cdot} = 2\{\cos ec^2\theta\cot\theta g(\tan\theta) - \sec^2\theta\tan\theta f(\tan\theta)\}\dot{\theta},$$
$$r^4\dot{\theta}^2 = L_\circ + 2\int\{\cos ec^2\theta\cot\theta g(\tan\theta) - \sec^2\theta\tan\theta f(\tan\theta)\}r^{-2}Ldt.$$

i.e. $$L^2 = L_\circ + \alpha(\theta) \tag{2.3}$$

where $L_\circ$ is a constant, and $L = r^2\dot{\theta}$ defined the angular momentum of the motion which is not constant in this dynamics.

Now by setting $u = r^{-1}$; $\dot{r} = -Lu_\theta$; $\ddot{r} = -L^2u^2u_{\theta\theta}$ and substituting into (2.1) we have that (Arunaye [11, 14, 22])

$$u_{\theta\theta} + [1 + \cos ec\theta h(\cot\theta) + \{\sec^2\theta f(\tan\theta) + \cos ec^2\theta g(\tan\theta)\}L^{-2}]u = 0.$$

i.e. $$u_{\theta\theta} + \omega^2(\theta)u = 0. \tag{2.4}$$

On taking $L_\circ = u_2$ and $u = u_1$ and imposing the conditions
$$f(\tan\theta) = \sin\theta\cos\theta = g(\tan\theta) = L^{-1} \tag{2.5}$$
such that $x = r\cos(\frac{1}{2}\sin^{-1}(2L^{-1}))$, $y = r\sin(\frac{1}{2}\sin^{-1}(2L^{-1}))$ and that $h(\cot\theta) \equiv 0$, when $H$ is a kernel transformation on the radial component of the motion. Then (2.4) and (2.3) respectively become

$$u_{1,\theta\theta} + 2u_1 = 0, \tag{2.6}$$
$$u_{2,\theta} = 0,$$

this is the reduced system for the Kepler-Ermakov systems.

The radial and transversal components of the motion for the generalized Ermakov systems (1.2) are given by

$$\ddot{r} - r\dot{\theta}^2 = \frac{1}{r^3}\{\sec^2\theta f(\tan\theta) + \cos ec^2\theta g(\tan\theta)\}, \tag{2.7}$$

$$r\ddot{\theta} + 2\dot{r}\dot{\theta} = -\frac{1}{r^3}\{\sec^2\theta\tan\theta f(\tan\theta) - \cos ec^2\theta\cot\theta g(\tan\theta)\}. \tag{2.8}$$

Similarly, we obtained the reduced system for the generalized Ermakov systems given by

$$u_{1,\theta\theta} + 2u_1 = 0, \tag{2.9}$$
$$u_{2,\theta} = 0$$

where the condition (2.5) is also assumed, with $L_\circ$ and $\alpha(\theta)$ defined by (2.3).

The radial and transversal components of the motion for the Ermakov-toy systems (1.3) are given by

$$\ddot{r} - r\dot{\theta}^2 = \frac{1}{r^3}(\tan\theta + \cot\theta)^2, \tag{2.10}$$

$$r\ddot{\theta} + 2\dot{r}\dot{\theta} = -\frac{1}{2r^3}(\tan\theta - \cot\theta)' \tag{2.11}$$

where prime implies derivation with respect to $\theta$. Similarly the same reduction procedure above produced

$$r^4\dot{\theta}^2 = L_\circ - \int\{\sec^2\theta - \cosec^2\theta\}r^{-2}L\,dt.$$

i.e. $\quad L^2 = L_\circ + \alpha(\theta).\tag{2.12}$

And on imposing the condition that the magnitude of the angular momentum satisfies $L = \tan\theta + \cot\theta$ with $2\theta = \sin^{-1}(2L^{-1})$ one obtains the following reduced system

$$u_{1,\theta\theta} + 2u_1 = 0, \tag{2.13}$$
$$u_{2,\theta} = 0$$

where $L_\circ = u_2$ and $u = u_1$, for the Ermakov-toy systems (1.3).

### 3.1 Lie point symmetries of the reduced three classes of Ermakov systems

We shall utilized (2.13) as hypothetical illustrative example of these classes of Ermakov systems under investigation to present the nonlocal symmetries of the Ermakov-toy systems and note that it is easy to deduce the nonlocal symmetries of the Kepler-Ermakov systems and generalized Ermakov systems from the following. Now we note that (2.13) is a system of second order linear equation and a conservation law. The system (2.13) has nine Lie point symmetries (well known in the literature). They are

$$\Gamma_1 = 2u_1\partial_1 + u_2\partial_2 \; ; \; \Gamma_2 = \partial_\theta \; ; \; \Gamma_3 = u_1\partial_1 \; ; \; \Gamma_{4\pm} = e^{\pm\sqrt{2}i\theta}\partial_1;$$
$$\Gamma_{6\pm} = e^{\pm 2\sqrt{2}i\theta}[\partial_\theta \pm iu_1\partial_1]; \; \Gamma_{8\pm} = e^{\pm\sqrt{2}i\theta}[u_1\partial_\theta \pm iu_1^2\partial_1]. \tag{3.1}$$

### 3.2 Generalized symmetries of the three classes of Ermakov systems

Substituting back for the original variables in the symmetries (3.1) of the Ermakov-toy systems we obtained the following generalized symmetries:

$$V_1 = [L^2 + \int \{\sec^2\theta - \csc^2\theta\} r^{-2} L dt] \partial_t - 2r^{-3} \partial_r ; \quad V_2 = r^{-2} L \partial_t ;$$
$$V_3 = -r^{-3} \partial_r ; \quad V_{4\pm} = -e^{\pm\sqrt{2}i\theta} r^{-2} \partial_r ; \quad V_{6\pm} = e^{\pm 2\sqrt{2}i\theta} [r^{-2} L \partial_t \mp ir^{-3} \partial_r] ;$$
$$V_{8\pm} = e^{\pm\sqrt{2}i\theta} [r^{-2} L \partial_t \mp ir^{-4} \partial_r] ; \quad V_{10} = \partial_t$$

where $L^2 = L_\circ + \alpha(\theta)$ as in (2.12) is obtained from (2.11) similarly as in (2.3) and $\theta = \int r^{-2} L dt$. The symmetry $V_{10}$ is introduced as the symmetry responsible for the reduction of order by the change of independent variable from $t$ to $\theta$.

The following generalized symmetries for the generalized Ermakov systems are obtained:

$$V_1 = [L^2 + 2\int \{\sec^2\theta \tan\theta f(\tan\theta) - \csc^2\theta \cot\theta g(\tan\theta)\} r^{-2} L dt] \partial_t - 2r^{-3} \partial_r ;$$
$$V_2 = r^{-2} L \partial_t ; \quad V_3 = -r^{-3} \partial_r ; \quad V_{4\pm} = -e^{\pm\sqrt{2}i\theta} r^{-2} \partial_r ;$$
$$V_{6\pm} = e^{\pm 2\sqrt{2}i\theta} [r^{-2} L \partial_t \mp ir^{-3} \partial_r] ;$$
$$V_{8\pm} = e^{\pm\sqrt{2}i\theta} [r^{-2} L \partial_t \mp ir^{-4} \partial_r] ; \quad V_{10} = \partial_t$$

where $L^2 = L_\circ + \alpha(\theta)$ as in (2.12) is obtained from (2.8) similarly as in (2.3) and $\theta = \int r^{-2} L dt$. The symmetry $V_{10}$ is introduced as the symmetry responsible for the reduction of order by the change of independent variable from $t$ to $\theta$.

We observe that these generalized symmetries are absolutely different from those obtained by Leach and Karasu (Alkanli) [7], Karasu (Alkanli) and Yildirim [16] and Leach, Karasu (Alkanli), Nucci and Andriopoulos [17] and Arunaye [22], consequences of their reduced systems. We also note that the imposed conditions on the magnitude of the angular momentum, with a specific function $H$ defining Kernel transformation on the radial component of motion influence the reduction of the Kepler-Ermakov systems to the generalized Ermakov systems.

### 4.1 Conclusions

In the forgoing, the reduction of three classes of the Ermakov systems to systems of two equations - one second order linear equation and an equation of motion was shown (Arunaye [14, 22]) after assuming certain conditions on the magnitude of the angular momentum of the motion. We note the variation in the reduced system (2.6), (2.9) and (2.13) from the usual reduced system of oscillator and conservation law. The application of Lie point symmetry algorithm produced the same nine point symmetries for each set of reduced systems however the backward transformation to original variables produced new generalized symmetries. Although the generalized symmetries obtained in this case may not be identical to those from Nucci reduction process, the variety is a consequence of the fact that nonlocal symmetries are infinite and have no unique algorithm for their general determination. We note (Arunaye [22]) the improvement on the number of symmetries from five obtained in Leach, Karasu (Alkanli), Nucci and Andriopoulos [17] to ten in this reduction process. We also observe that the conditions on the magnitude of the angular momentum of these dynamics may have some geometric implications on the dynamics, which may also affect the algebra of the complete symmetry group of these classes of Ermakov systems. This and exact symmetry transformations of these classes of Ermakov systems are subject for further discussion.